\RequirePackage{fix-cm}
\documentclass[smallextended]{svjour3}       
\smartqed  
\usepackage{amsmath}
\usepackage{amssymb}
\usepackage{graphicx}
\usepackage{amsfonts}
\usepackage[noend]{algpseudocode}
\begin{document}

\title{Comment on ``Approximation algorithms for quadratic programming"
\thanks{This research was supported by
	the Beijing Natural Science Foundation under grant Z180005, and
	the National Natural Science Foundation of China under grants 12171021 and 11822103.}
        }


\author{Tongli Zhang \and Yong Xia
}


\institute{
 T.L. Zhang \and Y. Xia \at
              LMIB of the Ministry of Education,
              School of Mathematical Sciences, Beihang University, Beijing,
100191, P. R. China
\email{yxia@buaa.edu.cn (Y. Xia, Corresponding author)}
}

\date{Received: date / Accepted: date}
\maketitle
\begin{abstract}
The radius of the outer Dikin ellipsoid of the
intersection of $m$ ellipsoids due to  Fu et al. (J. Comb. Optim., 2, 29-50, 1998) is corrected from $m$ to $\sqrt{m^2+m}$. The approximation bound for the general convex quadratic constrained nonconvex quadratic program is correspondingly corrected.
\end{abstract}


\section{Corrections}
Let $\mathcal{F}$ be the intersection of $m$ ellipsoids:
\begin{eqnarray*}
\mathcal{F}:=\left\{x\in\mathbb{R}^n:~g_i(x):=d_i-c_i^Tx-\frac{1}{2}x^TQ_ix\ge 0,~ i=1,\cdots,m\right\},
\end{eqnarray*}
where $Q_i\in\mathbb{R}^{n\times n}$ is symmetric and positive semidefinite, $c_i\in\mathbb{R}^n$ and $d_i\in\mathbb{R}$. The logarithmic barrier function for $\mathcal{F}$ is defined as
$
L(x):=-\sum\limits_{i=1}^m\log (g_i(x))
$ with gradient and Hessian being given by
\[
\nabla L(x)=\sum\limits_{i=1}^m\frac{c_i+Q_ix}{g_i(x)},~
\nabla^2L(x)=\sum\limits_{i=1}^m\left(\frac{\left(c_i+Q_ix\right)\left(c_i+Q_ix\right)^T}{g_i^2(x)}+\frac{Q_i}{g_i(x)}\right).
\]
The analytic center of $\mathcal{F}$ is defined as the unique minimizer of $L(x)$, denoted by $\bar{x}$.  The Dikin ellipsoid is defined as the following ellipsoid centered at $\bar{x}$ with radius $r>0$:
\begin{eqnarray*}
E\left(\bar{x};r\right)&:=&\left\{x\in\mathbb{R}^n:\left(x-\bar{x}\right)^T\nabla^2L\left(\bar{x}\right)\left(x-\bar{x}\right)\le r^2\right\}.
\end{eqnarray*}
It was proved in \cite{Nesterov1993} that
\begin{equation*}
E(\bar{x};1)\subseteq\mathcal{F}. 
\end{equation*}
As presented in Theorem $3$ in \cite{Fu1998}, Fu et al.  claimed that
\begin{equation}
\mathcal{F}\subseteq E(\bar{x};m). \label{m:2}
\end{equation}
It serves as a fundamental theorem in developing approximation algorithms \cite{Fu1998} for the convex quadratic constrained quadratic programming problem:
\begin{equation*}
{\rm(QP)}~~ {\rm min}_{x\in\mathcal{F}} ~ q(x):=\frac{1}{2}x^TQx+c^Tx.
\end{equation*}

With regret, a typo in the proof makes \eqref{m:2} incorrect.  It can be seen from the following simple example.
\begin{example}\label{exm1}
Let $\mathcal{F}_1:=\left\{x\in\mathbb{R}: x^2\le 1\right\}$.
The logarithmic barrier function is
$
L(x)=-\log(1-x^2)$ and $\nabla^2L(x)= 2(1+x^2)/(1-x^2)^2$.
It is not difficult to verify that the analytic center of $\mathcal{F}_1$ is $0$ and $\nabla^2L(0)=2$. Then we have
\[
E(0;1)=\{x\in\mathbb{R}^n: ~2x^2\le 1\}.
\]
Clearly, $E(0;1)\subseteq\mathcal{F}_1\nsubseteq E(0;m)=E(0;1)$.
\end{example}

Fortunately, \eqref{m:2} can be corrected as
\begin{theorem}\label{thm}
It holds that
\[
\mathcal{F}\subseteq E(\bar{x};\sqrt{m^2+m}).\]
\end{theorem}
The proof is very similar to that of Theorem $3$ in \cite{Fu1998}. For completeness, we leave it in the appendix.

Let us reconsider Example \ref{exm1}. We can verify that $\mathcal{F}_1=E(0;\sqrt{2})$, which implies that the corrected Theorem \ref{thm} is tight for this example.


For (QP) and $\epsilon\in[0,1]$, we call $x\in\mathcal{F}$ an $\epsilon$-minimizer if it satisfies that
\[
\frac{q(x)-\underline{z}}{\overline{z}-\underline{z}}\le \epsilon,
\]
where $\underline{z}$ and $\overline{z}$ are the minimum and maximum  of $q(x)$ over $\mathcal{F}$, respectively.

The approximation bound of Fu et al.'s approximation algorithm for solving (QP) (see Theorem 4 in \cite{Fu1998}) is correspondingly corrected as follows.

\begin{theorem}\label{thm2}
For all $\epsilon\in(0,1-1/\sqrt{2}]$ and all $m\ge 2$, there exists a polynomial-time approximation algorithm for computing an $(1-\frac{1-\epsilon}{c(m)(1+\epsilon)^2})$-minimizer of $q(x)$ over $\mathcal{F}$, where $c(m)=m^2+m$. 
\end{theorem}

When $m=1$, (QP) is known as the trust-region subproblem (TRS) that can be solved in polynomial time \cite{Fu1998}.
Celis-Dennis-Tapia (CDT) subproblem is a notable special case of (QP) with
$m=2$.  Theorem \ref{thm2} implies that an $\left(1-(1-\epsilon)/(6(1+\epsilon)^2)\right)$-minimizer for (CDT) can be computed in polynomial time.

\section{Appendix}
Proof of Theorem \ref{thm}. Denote $g_i(x)$, $g_i(\bar{x})$ and $\nabla L(\bar{x})$ by $g_i$, $\bar{g}_i$ and $\nabla \bar{L}$, respectively. Let
\begin{eqnarray*}
\Delta_i=\frac{1}{2}(x-\bar{x})^TQ_i(x-\bar{x}), ~i=1,\cdots,m.
\end{eqnarray*}
By Taylor expansion, we have
\begin{eqnarray}
-(c_i+Q_i\bar{x})^T(x-\bar{x})=g_i-\bar{g}_i+\Delta_i.\label{eq1}
\end{eqnarray}
Since $\bar{x}$ is the analytic center of $\mathcal{F}$, it follows from $\nabla L(\bar{x})=0$ and \eqref{eq1} that
\[
0=\sum\limits_{i=1}^m\frac{(c_i+Q_i\bar{x})^T(x-\bar{x})}{\bar{g}_i}=
m-\sum\limits_{i=1}^m\frac{g_i+\Delta_i}{\bar{g_i}},
\]
which further implies that
\begin{eqnarray*}
&&\sum\limits_{i=1}^m\frac{\Delta_i}{\bar{g}_i}=m-\sum\limits_{i=1}^m\frac{g_i}{\bar{g}_i}\le m,\\
&&m^2=\sum\limits_{i=1}^m\left(\frac{g_i+\Delta_i}{\bar{g}_i}\right)^2+2\sum\limits_{i\neq j}\left(\frac{g_i+\Delta_i}{\bar{g}_i}\right)\left(\frac{g_j+\Delta_j}{\bar{g}_j}\right)
\ge \sum\limits_{i=1}^m\left(\frac{g_i+\Delta_i}{\bar{g}_i}\right)^2.
\end{eqnarray*}
Therefore, we have
\begin{eqnarray}
(x-\bar{x})^T\nabla^2L(\bar{x})(x-\bar{x})&=&\sum\limits_{i=1}^m\left(\frac{((c_i+Q_i \bar{x})^T(x-\bar{x}))^2}{g_i^2(\bar{x})}+\frac{(x-\bar{x})^TQ_i(x-\bar{x})}{g_i(\bar{x})}\right)\notag\\
&=&\sum\limits_{i=1}^m\left(\frac{g_i+\Delta_i}{\bar{g}_i}-1\right)^2+\sum\limits_{i=1}^m\frac{2\Delta_i}{\bar{g}_i}\label{eq2}\\
&\le&\sum\limits_{i=1}^m\left\{\left(\frac{g_i+\Delta_i}{\bar{g}_i}\right)^2-
2\sum\limits_{i=1}^m\left(\frac{g_i+\Delta_i}{\bar{g}_i}\right)+1\right\}+2m\notag\\
&\le& m^2-2m+m+2m=m^2+m.\notag
\end{eqnarray}
The proof is completed. \hfill \quad{$\Box$}\smallskip

We remark that the authors in \cite{Fu1998} mistook $\sum\limits_{i=1}^m 2\Delta_i/\bar{g}_i$  for $\sum\limits_{i=1}^m \Delta_i/\bar{g}_i$ in (\ref{eq2}), which leads to the incorrect result.


\end{document}